\documentclass[12pt]{article}
\newcommand{\sss}{\setcounter{equation}{0}}

\newtheorem{theorem}{THEOREM}[section]

\newtheorem{remark}[theorem]{REMARK}

\newtheorem{prop}[theorem]{PROPOSITION}

\def\ER{{\bf R}}
\def\beq{\begin{equation}}
\def\ene{\end{equation}}
\def\bull{\begin{flushright} \vrule height 6pt width 6pt depth -.pt
\end{flushright}}

\def\0{W_{1,2}^{(0)}}
\def\1{\mathcal{H}_{1}^{(0)}}
\def\2{\mathcal{H}_{1}}
\def\3{\mathcal{H}_{2}}
\def\4{\mathcal{H}_2^{(0)}}
\def\u{e^{-itH}}
\def\l{\mathit{Lip}}
\begin{document}
%\baselineskip=21.6pt
% double space 28.8pt
%\baselineskip=28.8pt
\baselineskip=23.6pt
\title{The Forced Non-Linear Schr\"{o}dinger Equation with a Potential on the Half-Line
 \thanks{2000 {\sc AMS} classification 35Q40, 35Q55, and 35R30. Research partially supported
by proyecto PAPIIT, IN 101902-3, DGAPA-UNAM.}}
\author{ Ricardo Weder \thanks{Fellow, Sistema Nacional de Investigadores.}\\  Instituto de Investigaciones
en Matem\'aticas Aplicadas y en Sistemas\\ Universidad Nacional Aut\'onoma de M\'exico\\ 
Apartado Postal 20-726,
M\'exico D.F. 01000.      \\
E-Mail: weder@servidor.unam.mx\\}
\date{}
\maketitle

\begin{center}
\begin{minipage}{5.75in}
\centerline{{\bf Abstract}}\bigskip

In this paper  we prove that  the initial-boundary value  problem for the forced non-linear Schr\"{o}dinger 
equation with a potential on the half- -line is locally  and  (under stronger conditions) globally well posed, i.e. that there is a unique solution that depends continuously on 
the force at the boundary and on the initial data.  We allow for a large class of unbounded potentials. 
Actually, for local solutions  we have no restriction  on the grow at infinity of the positive part of the potential, 
and for global solutions very mild assumptions that allow, for example, for exponential grow.
 
\end{minipage}
\end{center}
\newpage
\section{Introduction}\sss
 In this paper we analyse in detail the 
initial-boundary value problem for the forced non-linear Schr\"{o}dinger equation with a 
potential  on the half-line (FNLSP),
\beq
i\frac{\partial}{\partial t}u(x,t)= -\frac{d^2}{dx^2}u(x,t)+ V(x)u(x,t)+ F(x,t,u),\, 
u(0,t)= f(t),\,
u(x,0)= \phi(x), 
\label{1.1}
\ene
where $F(x,t,u)$ is a complex-valued function of $x\in \ER^+, \, t\in 
\ER$, $u \in {\bf C}$. The functions $\phi$, $f$, satisfy the compatibility condition, $\phi(0)= f(0)$.
We solve this problem along the lines of \cite{ka}, who studied the pure initial value problem 
for the non-linear Schr\"{o}dinger equation on $\ER ^n, n \geq 1$.  Note, however, that we allow for a much larger class of potentials
than in \cite{ka}. In particular,  we do not need to require 
that $\frac{d^2}{dx^2}V(x)$  is bounded, as is the case in \cite{ka}. In fact, for local solutions  we have no restriction  on the grow at infinity of the positive part of the potential
and for global solutions very mild assumptions that allow, for example, for exponential grow.
 We consider potentials -that are in general time dependent- that can be decomposed as the sum of two parts. The first
one is what we call $V$ in  the FNLSP (\ref{1.1}); it is independent of time but it can  
have singularities and it can grow at infinity. The  second part is in general time dependent, but -together with its derivatives with 
respect to $x$ and $t$-  it has to be bounded 
for $(x,t) \in (\ER \times I)$, with $I$ any bounded set. This second part is included in $F$.

We consider the following class of potentials,
\beq
 V:= V_1 + V_2 \;  \mathrm{with} \; V_ j \in    L^1_{\mathrm{loc}}(\ER^+), j=1,2, V_1 \geq 0,  \mathrm{and }\sup_{x \in \ER^+} \int_x^{x+1} |V_2(y)|   \, dy < \infty.
\label{1.2}
\ene 
By $W_{l,2}, l=0,1,2,\cdots $, we denote the standard Sobolev
spaces \cite{ad} in $\ER^+$ and by  $ W_{l,2}^{(0)}$  the completion of 
$C^{\infty}_0(\ER^+)$ in the 
norm of $ W_{l,2}$. The functions in $ W_{l,2}^{(0)}$ satisfy the homogeneous 
Dirichlet boundary condition at zero, $\frac{d^j}{dx^j}u(0)=0, j=0,1,\cdots, l-1 $. In the case 
$l=0$ we use the standard notation, $W_{0,2}= W_{0,2}^{(0)}= L^2$. We designate, $q:= \sqrt{V_1}$, and by  $D(q)$ the domain  in $L^2$ of the operator 
of multiplication by $q$. We denote,
\beq
\mathcal{H}_{1}^{(0)} := \0 \cap D(q)     \; \mathrm{with \, norm} \; \|  \phi \|_{\1}:= \max \left[  \| \phi  \|_{\0} , \|  q \phi \|_{L^2} \right].
\label{1.3}
\ene  

Let us denote by $H_0$ the self-adjoint realization of $-\frac{d^2}{dx^2}$ with 
domain $ W_{2,2} \cap W_{1,2}^{(0)}$, i.e., the self-adjoint realization with homogeneous 
Dirichlet boundary condition at zero.  We have that (see Section 2 for details) the quadratic form,
\beq
\mathit{h}(\phi, \psi):= (\acute{\psi}, \acute{\psi})+ (V\phi, \psi), \mathrm{with \, domain}, \; D(\mathit{h}):= \1,
\label{1.4}
\ene
is closed and bounded below. We denote by $H$ the associated bounded-below, self-adjoint operator (see \cite{rs},  \cite{ka2}). Then,
$D(H) =\{ \phi \in  \1 : H_0 \phi + V \phi \in L^2\}  $ and,
\beq
H \phi = H_0 \phi + V \phi, \mathrm{for} \; \phi \in D(H).
\label{1.5}
\ene
%Let us take $M >0$  such that $H+M >0$. Then,
%\beq
%D\left( \sqrt{H+M} \right) = \1.
%\label{1.6}
%\ene 

We designate,
\beq
\mathcal{H}_{1} := W_{1,2} \cap D(q)     \; \mathrm{with \, norm} \; \|  \phi \|_{\2}:= \max \left[  \| \phi  \|_{W_{1,2}} , \|  q \phi \|_{L^2} \right],
\label{1.7}
\ene  
and,
$$
\mathcal{H}_{2} := \{ \phi \in \2 \;  \mathrm{such \; that}    (-\frac{d^2}{dx^2}+ V ) \phi \in L^2 \},
$$
\beq
 \mathrm{with \, norm} \; 
\|  \phi \|_{\3}:= \mathrm{max} \left[ \| \phi  \|_{\2} , \left\|\left(-\frac{d^2}{dx^2}+ V \right) \phi \right \|_{L^2}\right].
\label{1.8}
\ene

In Section 2 we study the initial-boundary value  problem for the FNLSP (\ref{1.1}). We first construct 
local solutions  assuming that for each fixed $x,t$, the non-linearity $F(x,t,u)$  is 
$C^1$ in the real sense as a function of $u$. We prove that the FNLSP (\ref{1.1}) is locally 
well posed in
$  \2$ and in $ \3$ and that  there is  continuos dependence on the initial and 
boundary data. In other words, the FNLSP (\ref{1.1}) forms a dynamical system by generating a continuous local 
flow (see \cite{ka}). Then, we prove that if $F$ satisfies a sign condition and has a hamiltonian
structure the solutions exist for all times. Under these conditions the continuous local flows become
global continuous flows, and in this sense the spaces of initial data $ \2$ and 
$ \3$ are fundamental for the FNLSP (\ref{1.1}).  Note that if $V_1 \equiv 0$, $\2 = W_{1,2}$. We give in Section 2 sufficient conditions  on $V$ assuring  that
$\3  =W_{2,2} \cap D(V_1) $ and in particular if $V_1 \equiv 0$, that $\3= W_{2,2}$.   
 
The existence and uniqueness of global solutions  in $W_{2,2}$ to the FNLSP (\ref{1.1}) with  
$V\equiv 0$ and  $F= \lambda |u|^2 u$ was proven in \cite{cb}, and the continuous dependence
on the initial value and the boundary condition in \cite{b2}. For existence and uniqueness of global  solutions with $V \equiv 0$ and  $F= \lambda |u|^{p-1} u, 
\lambda >0, p > 3 $ see  \cite{b1}. These papers give references for the application of the FNLSP  (\ref{1.1}) to important physical problems. 
For the solution of the direct and inverse scattering problems for the FNLSP (\ref{1.1}) see \cite{we1} 
and \cite{we}.
The existence of global solutions in $\ER^n, n \geq 2$, with $V \equiv 0$ and
$F = \lambda |u|^{p-1} u, 1 < p < \infty, \lambda > 0$, was proven in \cite{bs}. For the integrable case  where (\ref{1.1}) can be studied with inverse scattering transform 
methods see
\cite{fo} and the references quoted there. For the Korteweg-De Vries equation in the half-line see  \cite{bo} and   \cite{ke}.
For general references in non-linear initial value problems see \cite{st4}, \cite{ra}, \cite{gi} and \cite{bour}.

\section{The Initial Boundary-Value Problem}\sss
We first prepare results that we need. The Propositon below is well known. We give the simple proof for the reader's convenience
\begin{prop}
Suppose that
\beq
\sup_{x \in \ER}\int_x^{x+1}|V_2(y)|\, dy < \infty.
\label{2.1}
\ene

Then, for any $\epsilon >0$ there is a constant, $K_{\epsilon}$, such that,
\beq
\int_0^{\infty} |V_2(x)|\,  |\phi(x)|^2\, dx \leq \epsilon \left\|\acute {\phi} \,\right\|_{L^2}^2+ K_{\epsilon} \left\|\phi \right\|_{L^2}^2, \phi \in W_{1,2}.
\label{2.2}
\ene
Moreover, if 
\beq
\sup_{x \in \ER}\int_x^{x+1}|V_2(y)|^2\, dy < \infty,
\label{2.3}
\ene
 for any $\epsilon >0$ there is a constant, $K_{\epsilon}$, such that,
\beq
\left\| V_2 \phi \right\|^2_{L^2} \leq \epsilon \left\| H_0 \phi \right\|_{L^2}^2+ K_{\epsilon} \left\|\phi \right\|_{L^2}^2, \phi \in D(H_0).
\label{2.4}
\ene
\end{prop}

\noindent{\it Proof:} \,\,  If  $\phi \in W_{1,2 }$,  for any $n=0,1, \cdots$,
 any $x,y \in [n, n+1]$ and any $\delta > 0$, we have that,
\beq
|\phi(x)|^2- |\phi(y)|^2 = 2  \, \hbox{Re} \int_y^x \phi(z) \overline{\acute{\phi}(z)}\, dz
\leq \delta \int_n^{n+1}|\acute{\phi(z)}|^2 dz + \frac{1}{\delta} \int_n^{n+1}|\phi(z)|^2 dz.
\label{2.5}
\ene
By the mean value theorem we can choose $y$ such that,  $|\phi(y)|^2= \int_n^{n+1}|\phi(z)|^2 dz$, and it follows that,
\beq
|\phi(x)|^2 
\leq \delta \int_n^{n+1}|\acute{\phi}|^2(z) dz + \left(1+ \frac{1}{\delta}\right) \int_n^{n+1}|\phi(z)|^2 dz.
\label{2.6}
\ene
Let $C$ be the finite quantity in the left-hand side of (\ref{2.1}). Then,
\beq
\int_n^{n+1} |V_2(x)| |\phi(x)|^2\, dx \leq  C \delta  \int_n^{n+1}|\acute{\phi}|^2(z) dz + C \left(1+ \frac{1}{\delta}\right) \int_n^{n+1}|\phi(z)|^2 dz.
\label{2.7}
\ene
Taking $\delta$ so small that  $  \epsilon= \delta C$, and adding over $n$  we obtain (\ref{2.2}). 
Let us now denote by $C$ the finite quantity on the 
left-hand side of (\ref{2.3}). As $D(H_0)=  W_{2,2} \cap \0$, if follows from  (\ref{2.6}) that, 
\beq
\int_n^{n+1} |V_2(x)|^2 |\phi(x)|^2\, dx \leq  C \delta  \int_n^{n+1}|\acute{\phi}|^2(z) dz + C \left(1+ \frac{1}{\delta}\right) \int_n^{n+1}|\phi(z)|^2 dz.
\label{2.8}
\ene
Taking now $\delta$ so small that  $  \epsilon= \delta C/2$, adding over $n$ and as $ \| \acute{\phi}\|_{L^2}^2= (H_0 \phi, \phi)
\leq \|H_0 \phi\|_{L^2}^2 /2 + \| \phi\|_{L^2}^2 /2$,  we obtain (\ref{2.4}).  

\bull
Assuming that (\ref{1.2}) holds, the results about  $\mathit{h}$ and $H$ stated in the introduction (see (\ref{1.4})-(\ref{1.5}))   follow from (\ref{2.2}) and 
 \cite{rs}, \cite{ka2}. Below we always assume that (\ref{2.1}) is satisfied.

We  study the initial boundary-value problem for the FNLSP (1.1) for $ t \geq 0$, but by changing $t$ into
$-t$ and taking the complex conjugate of the solution (time reversal) we also obtain the results
for $t \leq 0$.  
Let $F(x,t,z)$ be a complex-valued function of $ x\in \ER^+,\, t \in \ER^+,\, z\in {\bf C}$. As 
we are not assuming analyticity of $F$ we consider the derivative, $\acute{F}$, in the real sense.
 For each $z\in {\bf C}$, $\acute{F}$ is defined as the real-linear operator on ${\bf C}$, given by, 
\beq
\acute{F}(x,t,z)v:= \left(\frac{\partial}{\partial z}\,F(x,t,z)\right)v + 
\left(\frac{\partial}{\partial \overline{z}}F(x,t,z)\right)\overline{v}, \,\, v \in {\bf C},
\label{3.1}
\ene
 with the standard notation, $\frac{\partial}{\partial z}:=(1/2)\left(\frac{\partial}{\partial a}
-i \frac{\partial}{\partial b }\right)$ and  $\frac{\partial}{\partial \overline{z}}
:=(1/2)\left(\frac{\partial}{\partial a}+
i \frac{\partial}{\partial b }\right)$ where $z=a+ib$. We denote, $\left|  \acute{F}\right|:=
\left|\frac{\partial}{\partial z} F\right|+ \left|\frac{\partial}{\partial \overline{z}}
F\right|$.  $\acute{F}$ can be identified (when viewed as a $2 \times 2$ matrix) with the Gateaux derivative
in the real sense of the map $ z \in {\bf C} \rightarrow F(x,t, z)\in {\bf C}$ for each fixed
$x,\, t$. 
We  say that for each fixed $x,t$, $F$ is $C^1\left({\bf C}, {\bf C}\right)$  in the
real sense if $\frac{\partial}{\partial z}\,F(x,t,z)$ and 
$\frac{\partial}{\partial \overline{z}}F(x,t,z)$ are continuous functions of $z$ for each fixed
$x,t$, or equivalently if the map $z \rightarrow \acute{F}(x,t,z)$ is continuous from ${\bf C}$
into the real-linear operators in ${\bf C}$. 
For $T >0$   we denote, $I:=[0, T]$ if $  T < \infty$ and $I:=[0, \infty)$ if 
$T= \infty$.

\noindent {\bf Assumption A}

\noindent Suppose that $F(x,t,z)$ is a function from $ [0, \infty ) \times I \times
{\bf C}$ into ${\bf C}$, that for each fixed $x \in [0, \infty ), t \in I$, is $C^1$ in $z$ in the real sense. Moreover, assume that  for 
each fixed $t,z$, $F$ is differentiable in $x \in \ER^+$, that $F(x,t,0)=0$ and that for each $R >0$ 
and each bounded subset, $I_N$, of $I$, there is a constant $C_{R,N}$ such that,
\beq
\left|\acute{F}(x,t,z)\right|
\leq \, C_{R,N}, \,\,\hbox{for}\,\, x\in [0, \infty),\, t\in I_N, \, |z|\leq R,
\label{3.2}
\ene
and,
\beq
\left|\left(\frac{\partial}{\partial x}F\right)(x,t,z)\right|
\leq \, C_{R,N} \,\, |z|, \,\, \hbox{for}\,\, x\in [0, \infty),\, t\in I_N, \, |z|\leq R.
\label{3.2b}
\ene
Furthermore, if the force, $f$, in (\ref{1.1}) is not identically zero, suppose that for each fixed $ z$, $ F(0,t,z)$
is differentiable in $t$ and
\beq
\left|\frac{\partial}{\partial t}F(0,t,z)\right|\leq C_{R,N}, \, t \in I_n, |z|\leq R.
\label{3.2c}
\ene 

\bull
Assumption A allows for a large class of non-linearities. For example, the standard single-power 
non-linearity, $|u|^{p-1}\,u, p > 1$, or more generally any $F(z)$ that is $C^1$ in the 
real sense, are allowed. 
Let us denote by $ {\mathcal H}_{-1}$ the dual of $\1$ with the pairing given by the scalar 
product of $L^2$. Then, as the quadratic form domain of $H$ is    $ \1$,
 $H$ extends to a bounded operator from  $\1$ into $ {\mathcal H}_{-1}$.
Moreover, $\u$ is a bounded operator from $\1$ into 
$C\left( I, \1\right) \cap C^1\left( I, {\mathcal H }_{-1}\right)$ and,
\beq
i\frac{\partial}{\partial t} \u \, \phi = H \u\phi= \u H \phi.
\label{3.6}
\ene

Suppose that $u(x,t)\in C\left( I, \2  \right)$ is a solution of (\ref{1.1}) where
$f \in C^{2}(I)$. 
Furthermore, if $f$ is not identically zero, assume  that $ V \in W_{1,2}((0,\delta ))$ for 
some $ \delta >0$.
Note that the compatibility condition $\phi(0)= f(0)$ has to be satisfied if there is a 
solution to (\ref{1.1}).
Denote, $v(x,t):= u(x,t)- r(x,t)$ where $r(x,t):=[f(t)+ \frac{1}{2}x^2(V(0)f(t)+ F(0,t,f(t))-
i \acute{f}(t))] g(x)$, with $g \in C^{\infty}_0([0, \infty))$,
 $g(x)=1, 0 \leq x \leq \delta/2 $ and with support contained in $[0, \delta )$. Then, 
$v(x,t) \in 
C\left( I, \1  \right)$ solves,
\beq
i\frac{\partial}{\partial t}v(x,t)= Hv(x,t) + F_1(x,t,v),\, 
v(0,t)= 0,\, v(x,0)= v_0(x):=\phi(x)- r(x,0),
\label{3.7}
\ene
 
where,
\beq
F_1(x,t, v):= F(x,t, v+r) - i \frac{\partial}{\partial t}r + V(x)\, r- 
\frac{\partial^2}{\partial x^2}r.
\label{3.8}
\ene
Note that by the compatibility condition, $v_0 \in \1$. Clearly, equations (\ref{1.1}) and 
(\ref{3.7})
are equivalent. By Assumption A and Sobolev's \cite{ad}
theorem for any $t, N >0$ there is a constant, $C$, such that, 
$\left\|F_1(x,s,v ) \right\|_{L^2}\leq
C $ for all $0 \leq s \leq t$ and $v \in W_{1,2}$ with $\left\| v \right\|_{W_{1,2}}
\leq N$. Multiplying both sides of (\ref{3.6}) (evaluated at $\tau$) by $e^{-i(t-\tau)H }$ and integrating in $\tau$ from zero
to $t$ we obtain that,
\beq
v(t)= \u v_0 + \frac{1}{i} G {\cal F}(v),
\label{3.9}
\ene
where we designate  by ${\cal F}$ the operator 
$ v \rightarrow {\cal F}(v):= F_1(x,t, v)$, and
\beq
\left( Gf \right)(t):= \int_0^{t}\, e^{-i(t-\tau)H} \, f(\tau)\, d\tau,
 \, t \in I.
\label{3.9b}
\ene   
We prove below that if $v(t)\in C\left( I, \1  \right)$ is a solution to (\ref{3.9}),
it is also a solution to (\ref{3.7}). 
We denote,  
\beq
v_1(x,t):= \frac{1}{i}\int_0^t \, e^{-i (t- \tau )H} F_1(x,\tau, v(x,\tau))
\, d\tau.                        
\label{3.11}
\ene 
It follows from Assumption A that $ F_1(x,\tau, v) 
\in L^{\inf}_{\mathrm loc}\left( I, \1\right)$. Hence, $v_1 \in  
C\left( I, \1\right) \cap C^1\left( I, {\mathcal H}_{-1}\right)$, and
\beq
i\frac{\partial}{\partial t}v_1(t)=  H v_1(t)+  F_1(x,t, v).
\label{3.13}
\ene
Equations (\ref{3.9}) and (\ref{3.13}) imply that (\ref{3.7}) holds.
This proves that (\ref{3.7}) and (\ref{3.9}) are equivalent. We obtain our results
below by solving the integral equation (\ref{3.9}). 

\noindent {\bf Assumption B}

\noindent
Suppose that $V$ can be decomposed as,
$ V:= V_1 + V_2$  with   $ V_ j \in    L^1_{\mathrm{loc}}(\ER^+), j=1,2,  \;  V_1 \geq 0$, and $V_2$ satisfies (\ref{2.1}). 
Moreover, assume  that $f \in C^{2}(I)$, 
and  if $f$ is not identically zero, suppose  that $ V \in W_{1,2}((0,\delta))$ for some
$\delta >0$.

\bull

\noindent
We designate,\, ${\cal M}:= L^{\infty}\left(I, \1\right)$. By Sobolev's theorem \cite{ad},
${\cal M} \subset L^{\infty}\left(  \ER^+  \times I \right)$  and 
$\|v(x,t)\|_{L^{\infty}(\ER^+ \times I)}\leq C \|v\|_{{\cal M}}$. Moreover, we denote, $B:= L^{\infty}(I, L^2)$ and $\acute{B}:= L^1(I, L^2)$.

\begin{theorem}
Suppose that  Assumptions A and B  are satisfied. Then, for any $\phi \in \2$ satisfying
$\phi(0)= f(0)$,  there is a 
finite $T_0 \leq
T$ such that the FNLSP (\ref{1.1}) has a unique solution, $u \in C\left([0,T_0 ], 
\2\right)$
with $u(x,0)= \phi$. $T_0$ depends only on $\| \phi\|_{\2}$.
\end{theorem}

\noindent {\it Proof:}\,\, we prove the theorem by showing that (\ref{3.9}) has a unique 
solution $v \in C\left([0,T_0 ], \1\right)$ such that, $v(x,0)= v_0(x):=\phi(x)- r(x,0)$. Let us take $I=[0,T_0]$ with $T_0 < \infty$.
Let us denote by $\overline{{\cal M}}$ the space of bounded and continuous functions from 
$I$ into $\1$. 
Let ${\cal M}_R$ and $\overline{{\cal M}}_R$ be, respectively, the closed ball in ${\cal M}$,
and in $\overline{{\cal M}}$,  with center zero and radius $R$. Let us prove that ${\cal M}_R$
is closed in the norm of $B$. Suppose that $ v_n \in {\cal M}_R$ converges to $v$ in the 
norm of $B $. Then, $\lim_{n \rightarrow \infty}\|v_n(t)-v(t)\|_{L^2}=0$ for a.e. $t$.
But as $v_n \in {\cal M}_R$, $\|v_n(t)\|_{\1}   \leq R$ for a.e.
$t$. In consequence $\|v(t)\|_{\2}  \leq R$ for a.e.
$t$. Moreover, there is a subsequence - denoted also $v_n(t)$-  
 such $\frac{\partial}{\partial x}v_n$  converges weakly to 
$\frac{\partial}{\partial x}v$ in $L^2$, for a.e. $t$, and then, $ v(x,t)= 
\lim_{n \rightarrow \infty}
v_n(x,t)= \lim_{n \rightarrow \infty} \int_0^{x}\frac{\partial}{\partial y}{v}_n(y,t) dy = 
\int_0^{x}\frac{\partial}{\partial y}v(y,t)dy$,
and it follows that $v(0,t)=0$, i.e., $v \in {\cal M}_R$. Hence, ${\cal M}_R$ is a complete 
metric space in the norm of
$B$.
  
We define,
\beq
P(v):= \u v_0 + \frac{1}{i} G {\cal F}(v), v \in {\cal M}.
\label{3.16}
\ene
As $D(\sqrt{H+M}) = \1$,  for $M$ large enough, we have that the norm $\| \sqrt{H+M} \phi  \|_{L^2}$ is equivalent to the norm of $\1$. Then, by the unitarity of $\u$ in $L^2$ 
and as $ \sqrt{H+M}$ commutes with $\u$,   
\beq
\|\u v_0\|_{\1} \leq C \|\sqrt{H+M}\, \u v_0\|_{L^2} \leq C \|v_0 \|_{\1},
\label{3.17}
\ene 
and moreover,
\beq
\left\| \u \right\|_{{\mathcal B}\left(\1\right)} \leq C, \;\mathrm {for } \;t \in \ER.
\label{3.17b}
\ene
As $V\in W_{1,2}((0,\delta))$ if $f$ is not identically zero, 
it follows from Assumption A and Sobolev's theorem that there is a constant $C_R$
such that for $ v \in {\cal M}_R$,
\beq
\left\|{\cal F}(v) \right\|_{\1}\leq \ C_R.
\label{3.19}
\ene
Note that $ {\cal F}(v)(t)\in \1$ if
$ v \in {\cal M}$. By  (\ref{3.17b}) and  (\ref{3.19}), there is a constant $C_R$ such that,

\beq
\left\| G{\cal F}(v) \right\|_{{\cal M}}=\sup_{t \in [0,T_0]}  \left\|G{\cal F}(v)
\right\|_{\1}\leq  C_R \int_0^{T_0} \,dt \, \left\|{\cal F}(v)\right\|_{\1}\leq C_R \, T_0,
\label{3.20}
\ene  
for all $v \in {\cal M}_R$. By (\ref{3.17}) and (\ref{3.20}) we can take  $R$ large
enough and  $T_0$ small enough (depending only on $\|\phi \|_{\1}$) such that  $P$ sends 
${\cal M}_R$ into $\overline{{\cal M}_R}$. By Assumption A there is a constant $C$ such that,
\beq
\left\|{\cal F}(u)-{\cal F}(v)\right\|_{B}\leq C \left\|u-v  \right\|_{B},\, u,v\in 
{\cal M}_R.
\label{3.21}
\ene
Then, by the unitarity of $e^{-itH}$ in $L^2$,
\beq
\left\| P(u)- P(v)\right\|_{B}\leq C\,T_0 \, \left\| u-v  \right\|_{B},\, u,v\in 
{\cal M}_R.
\label{3.22}
\ene
Given $R$ we can take $T_0$ so small that $P$ is a contraction on the metric of $B$.
By the contraction mapping theorem $P(u)$ has a unique fixed point that is the only solution
to the FNLSP in ${\cal M}_R$ . If there is another solution, 
$u_1 \in C\left([0,T_0], \1\right)$, to (\ref{1.1}), then, $v_1:= u_1- r$ has to be a solution to 
(\ref{3.9}), but since $ \left\|G {\cal{F}}(v_1)(t)\right\|_{\1}$ can be made arbitrarily small 
by taking $ 0 \leq t \leq T_1, T_1 \leq T_0$, by the same argument as above we
have that $ v(t)= v_1(t), 0 \leq t
\leq T_1$, where $T_1$ depends only on $\|\phi\|_{\2}$. By iterating this argument we prove
that $v(t)=v_1(t), 0 \leq t \leq T_0$.

\begin{theorem}
Suppose that Assumptions A and B are satisfied by $ V, F, f, f_n, n=1,2, \cdots $, where we 
require that 
$ V \in W_{1,2}((0, \delta))$ for some $\delta >0$ only if the  $f_n, n=1,2, \cdots $ are not 
all identically zero for n large enough.
Then, the solution \, $u\in C\left([0,T_0], \2 \right), u(0)= \phi,\,  T_0 \leq T$, to the FNLSP (\ref{1.1}) 
depends continuously on the initial value and
on the boundary condition. In a precise way, let $u\in C\left([0, T_0], \2 \right)$ be the 
solution
to (\ref{1.1}) with $u(0)= \phi$, let $ \phi_n \rightarrow \phi$ in $\2$ with $ \phi_n(0)= f_n(0$) and 
assume that 
$ f_n(t) \rightarrow f(t)$ in 
$C^2([0,T_0])$.
Then, for $n$ large enough the solution $u_n \in C\left([0,T_0], \2 \right) $ to the 
FNLSP (\ref{1.1}) with initial condition
$\phi_n$ and boundary condition $f_n$ exits for $t \in [0,T_0]$ and $u_n \rightarrow u$ in 
$C\left([0,T_0], \2 \right)$.  
\end{theorem}

\noindent {\it Proof:}\,\, We first prove a local version for $T_0$ small enough. We denote
$v_{0,n}(x):= \phi_n(x) - r_n(x,0),$ with $r_n(x,t):= [f_n(t)+ \frac{1}{2}x^2(V(0)
f_n(t)+ F(0,t,f_n(t))-
i \acute{f_n}(t))] \; g(x)$,  and,
 \beq
P_n(v):= \u v_{0,n} + \frac{1}{i} G {\cal F}_n(v), v \in {\cal M},
\label{3.23}
\ene
where we designate  by ${\cal F}_n$ the operator 
$ v \rightarrow {\cal F}_n(v):= F_{1,n}(x,t, v(x,t))$, with,
\beq
F_{1,n}(x,t,v):= F(x,t, v+r_n)- i \frac{\partial}{\partial t}r_n(x,t) + V(x)\, r_n(x,t)- 
\frac{\partial^2}{\partial x^2}r_n(x,t).
\label{3.24}
\ene
As in the proof of Theorem 2.2 we prove that for $R$ large enough and $T_0$ small enough
all the $P_n$ send ${\cal M}_R$ into $\overline{{\cal M}_R}$ and are contractions in the
norm of $B$ with a uniform contraction rate $\sigma < 1$ independent of $n$. Let $v_n$
be the unique fixed point. Then, $ u_n:= v_n+ r_n \in C\left([0,T_0], \2 \right)$
are the unique solutions to the FNLSP (\ref{1.1}) with initial value $\phi_n$ and boundary 
condition $f_n$.
Furthermore, 
$$
\|v_n - v \|_{B}=\left\|e^{-it H}[v_{0,n} -v_0]  + 
\frac{1}{i} \left(G {\cal F}_n(v_n)- G {\cal F}(v)\right)\right\|_{B} \leq \|v_{0,n} 
-v_0 \|_{L^2}+
 \sigma \left\| v_n - v \right\|_{B}+
$$
\beq
 C \, T_0 \|f_n -f\|_{C^2([0,T_0])},
\label{3.25}
\ene                      
and it follows that $ v_n \rightarrow v$ in $B$. 
Moreover, by (\ref{3.17}), (\ref{3.17b}) and denoting, $v_x:= \frac{\partial}{\partial x}v$, 
and $v_{n,x}:= \frac{\partial}{\partial x} v_n$,
$$
\left\| v_n(t)-v(t)\right\|_{\1} \leq  C \|v_{0,n} -v_{0} \|_{\1}+
 C T_0 \left[\left\| {\cal F}_n(v_n)- {\cal F}(v) \right\|_{B}+
\left\|D_n(v_n, v_{n,x}) -D_n(v_n,v_x)\right\|_{B}+ \right.
$$
\beq 
\left.
\left\|D_n(v_n, v_x) -
D(v,v_x)
\right\|_{B} + \left\| q({\cal F}_n(v_n)- {\cal F}(v)) \right\|_{B}\right],
\label{3.26}
\ene
where,
$$
D(v,v_x):= \acute{F}(x,t, v+r)(v_x + \frac{\partial}{\partial x}r) + 
\left(\frac{\partial}{\partial x}F
\right)(x,t, v+r)
$$
\beq
 - i \frac{\partial^2}{\partial x \partial t}r + \acute{V}(x)\, r+ V(x)\, 
\frac{\partial}{\partial x}r- 
\frac{\partial^3}{\partial x^3}r, 
\label{3.27}
\ene
and,
$$
D_n(v_n,v_{n,x}):= \acute{F}(x,t, v_n+r_n)(v_{n,x}+ \frac{\partial}{\partial x}r_n) + 
\left(\frac{\partial}{\partial x}F
\right)(x,t, v_n+ r_n)- 
$$
\beq
  - i \frac{\partial^2}{\partial x \partial t}r_n + \acute{V}(x)\, r_n + V(x)\, 
\frac{\partial}{\partial x}r_n - 
\frac{\partial^3}{\partial x^3}r_n. 
\label{3.28}
\ene 
Furthermore,
\beq
\left\|{\cal F}_n(v_n) - {\cal F}(v)\right\|_{B}+\left\|q({\cal F}_n(v_n) - {\cal F}(v)) \right\|_{B}  \leq C [
 \left\| v_n - v \right\|_{B}+\left\| q( v_n - v ) \right\|_{B}+  \|f_n -f\|_{C^2([0,T_0])}].
\label{3.29}
\ene
Also,
\beq
\left\|D_n(v_n,v_{n,x} )- 
D_n(v_n,v_x) \right\|_{B}\leq C \,  \left\| v_{n,x}-v_x \right\|_{B}.
\label{3.30}
\ene
But then,
$$
\left\| v_n-v\right\|_{C\left([0,T_0], \1\right)} \leq C \|v_{0,n} -v_{0} \|_{\1}+ 
 C T_0\left[  \left\| v_n(t)-v(t)
\right\|_{C\left([0,T_0], \1\right)}+
     \|f_n -f\|_{C^2([0,T_0])}+ \right.
$$
\beq
 \left.
\left\| D_n(v_n, v_x) -
D(v,v_x)
\right\|_{B}\right]. 
\label{3.31}
\ene

And it follows that for $ C T_0 < 1/2 $,
$$
\lim_{n \rightarrow \infty}\left\| v_n-v_x \right\|_{C\left([0,T_0\right], \1)}
 \leq 2 \,\,
 C \lim_{n 
\rightarrow \infty} \left[\left\| v_{0,n}- v_{0}\right\|_{\1}+  T_0 \left(\,\left\|
D_n(v_n,v_x )- 
D(v,v_x) \right\|_{B}+
\right.\right.
$$
\beq
\left.
\left.
 \|f_n -f\|_{C^2([0,T_0])}\right) \right]=0,
\label{3.32}
\ene
where we used that, as $ v_n \rightarrow v$ in $L^2$ and $\| v_n \|_{\0}\leq C$, it follows by 
interpolation
\cite{rs} that $v_n \rightarrow v$ in $W_{s,2}^{(0)}, 0 < s <1$, and Sobolev's theorem.
This proves that $ u_n \rightarrow u$ in $  C\left([0,T_0], \2\right)$. In a standard way 
we 
extend the result of the theorem -step by step- to the original interval. For this purpose it 
is essential that the interval 
of existence given by Theorem 2.2 depends only on the $\2$ norm of $\phi$. 

\begin{remark}
{\rm
\noindent Suppose that Assumptions A and B are satisfied with $I=[0, \infty)$. Let $T_m$ be the 
maximal time such that the solution, $u$, given by Theorem 2.2 can be extended
to a solution  $u \in C\left( [0, T_m), \2 \right)$, to the FNLSP (\ref{1.1}) with $ u(0)= \phi$.
Then if $T_m$ is finite we necessarily have that $\lim_{t \uparrow T_m} 
\,\|u(t)\|_{\2}= \infty$. In other words, the solution exists for all times unless it blows up
in the $\2$ norm for some finite time. This result  follows from Theorem 2.2, because if 
$\|u(t)\|_{\2}$ remains 
bounded  as $ t \uparrow T_m$ we can extend the solution $u$ continuously to $T_m + \epsilon$
for some $\epsilon >0$, contradicting the definition of $T_m$.
Theorem 2.2  implies also that the FNLSP (\ref{1.1}) has a unique solution, 
$u \in C\left( I, \2 \right)$, with $u(0)= \phi$. For, suppose that there is another solution,
$v \in C\left( I, \2 \right)$, of this problem. Then, by Theorem 2.2 $u(t)=v(t)$, for 
$ t \in I_0:=[0,T_0], 0 < T_0 \leq T$. Let $I_m :=[0, T_m) \subset I$ 
be the maximal interval such that $u(t)=v(t), t \in I_m$. Then, if 
$ T < \infty, \, T_m=T$ -note that by continuity this implies that $ u(T)= v(T)$-
and if $T= \infty, \, T_m = \infty$. If $T < \infty$  this follows because by Theorem 2.2 if $T_m < \infty$,  
$u(t)=v(t)$ for $t \leq T_m+ \epsilon$, for some $\epsilon > 0$, contradicting the definition 
of $T_m$.  By the same argument if $T= \infty$, 
$T_m$ can not be finite.}
\end{remark}  

\bull

We now consider solutions in $\3$.

\begin{theorem}
Suppose that  Assumptions A  and B are  satisfied. Furthermore, assume that for each fixed $x,z$, $F(x,t,z)$ is 
differentiable in $t$ and,
\beq
\left|\left(\frac{\partial}{\partial t}F\right)(x,t,z)\right|
\leq \, C_{R,N} |z|, \,\,\hbox{for}\,\, x\in [0, \infty ),\, t\in I_N, \, |z|\leq R.
\label{3.35f}
\ene
Then, for any $\phi \in \3$ with $\phi(0)= f(0)$ there is a 
finite $T_0 \leq
T$ such that the FNLSP (\ref{1.1}) has a unique solution $u \in C\left([0,T_0 ], 
\3\right)$
with $u(x,0)= \phi$. $T_0$ depends only on $\| \phi\|_{\3}$.
\end{theorem}

\noindent {\it Proof:}\,\, We designate,
\beq
{\mathcal H}_2^{(0)}:= \{ \phi \in \3: \phi(0)=0  \}, \; \mathrm {and} \;                                  
{\cal N}:= \left\{ v \in L^{\infty}\left( [0, T_0], \4 \right): 
\frac{\partial}{\partial t}v (t)\in B
\right\},
\label{3.36}
\ene
with norm
\beq
 \|v\|_{{\cal N}}:= \mathrm{max}\left[\|v\|_{L^{\infty}\left( [0, T_0], \4\right) },
\left\|\frac{\partial}{\partial t}v\right\|_{B}\right].
\label{3.37}
\ene 
We define ${\overline{\cal N}}$ as in (\ref{3.36}) but replacing $L^{\infty}$ with continuous.
Note that if  $v \in B$ and  
$\frac{\partial}{\partial t}v (t)\in \acute{B}$,
it follows that $v(t)$ is a absolutely continuous function of $t \in [0, T_0]$, with values 
in $L^2$. In consequence, $ v(0)\in L^2$ exists and,
\beq
\|v(0)\|_{L^2} \leq \| v \|_{B}.
\label{3.38} 
\ene
We use the designation,
\beq
{\cal N}_R:=\left\{ v \in {\cal N}: \|v\|_{{\cal N}}\leq R, \, \hbox{and}\, v(0)= v_0 \right\}.
\label{3.39}
\ene
We first prove 
that ${\cal N}_R$ is a complete metric space in the norm of $B$. It is enough to prove that
it is a closed subset of $B$. Supose that $v_n \in {\cal N}_R$ converges to $v \in B$ in 
the norm
of $B$. We have to prove that $ v \in {\cal N}_R$. We have that 
$\lim_{n \rightarrow \infty}\|v_n(t)-v(t)\|_{L^2}=0,\,  \hbox{for a.e.}\, \,t$.
But as $v_n \in {\cal N}_R$, $\|v_n(t)\|_{\cal N} \leq R$ for a.e.
$t$. In consequence,  $  \mathrm{max}  [[\|v\|_{L^{\infty}\left( [0, T_0], \3\right) },
\left\|\frac{\partial}{\partial t}v\right\|_{B}]
\leq R$ for a.e. $t$. We prove that $ v(0,t)=0$
as in the proof of Theorem 2.2. 
Moreover, we have that (eventually passing to a subsequence) $ \frac{\partial}{\partial t}v_n
\rightarrow \frac{\partial}{\partial t} v$ weakly. Then, as $ v_n(t)= v_0 + \int_0^t
 \frac{\partial}{\partial s}v_n(s)\, ds$, we obtain that $v(0)= v_0$. Hence, $v \in{\cal N}_R$.

Let  $P$ be defined as in (\ref{3.16}).
Let us prove that we can take $R$ so large and $T_0$ so small (depending only on 
$\|\phi\|_{\3}$)
that $P$ sends ${\cal N}_R$ into $\overline{{\cal N}_R}$. As $D(H)= \4$, and since $ H$ commutes with $\u$ and  $ i\frac{\partial}{\partial t}\u \phi = H \u \phi$,
$\u$ is bounded from $\4$  into $ \overline{{\cal N}}$ with operator norm independent of $T_0$.
Furthermore, suppose that $w \in B$ and that  
$\frac{\partial}{\partial t}w (t)\in \acute{B}$, with $w(0)= \psi \in L^2$.
Then,
\beq
\frac{\partial}{\partial t} G w = G \frac{\partial}{\partial t}w + \u \psi.
\label{3.40}
\ene
We write, $G w = \u w_1(t)$, with $w_1(t):= \int_0^t e^{i\tau H }\, w(\tau)\, d\tau$.
Then,
\beq
 \frac{\partial}{\partial s}e^{is H}w_1(t)\big|_{s=0}=e^{itH}w(t)- \psi -\int_0^t e^{i\tau H}
\, \acute{w}(\tau) \, d\tau.
\label{3.41}
\ene
As the right-hand side of (\ref{3.41}) belongs to $L^2$ for a.e. $t$, it follows that  $ w_1(t)\in D(H)= \4$ for a.e. $t$.
Then, $G w = \u w_1(t) \in \4$ for a.e. $t$, and
\beq
H G w = i \frac{\partial}{\partial t} G w - iw.
\label{3.42}
\ene
  By (\ref{3.40}) and (\ref{3.42}) $ G w \in {\cal N}$ and,
\beq
\left\| G w    \right\|_{{\cal N}}\leq C \left[\|w\|_{B}+ \left\|
\frac{\partial}{\partial t} w 
\right\|_{\acute{B}}\right].
\label{3.43}
\ene 

For $v \in {\cal N}_R$ we write
\beq
P(v)= \u v_0- \frac{1}{i} G{\cal F}(v_0)+ \frac{1}{i}G [{\cal F}(v)- {\cal F}(v_0)].
\label{3.44}
\ene

We take $R$ so large and $T_0$ so small  that $\|\u v_0 - \frac{1}{i} G {\cal F}(v_0)
\|_{{\cal N}} \leq R /2 $. Here we take $w = {\mathcal F}(v_0)$ in the estimates above.
We now put $w(t):= {\mathcal F}(v)(t)- {\cal F}(v_0)$.  By Assumption A, $ w \in B$ and
$\frac{\partial}{\partial t}w \in \acute{B}$. Then, as $w(0)=0$, by (\ref{3.43}) given $R$ we can take
$T_0$ so small that,
\beq
\left\|\frac{1}{i}G [{\cal F}(v)- {\cal F}(v_0)]\right\|_{{\cal N}_R} \leq \frac{R}{2}.
\label{3.45}
\ene 
With this choice of $R$ and $T_0$, $P$ sends ${\cal N}_R$ into $\overline{{\cal N}}_R$. We already know 
-see the proof of Theorem 2.2- that $P$ is a contraction in the norm of $B$. The unique fixed
point is the only solution to the FNLSP (\ref{1.1}) in $\overline{{\cal N}}_R$. We complete the proof of the 
theorem as in the proof of Theorem 2.2.

\begin{theorem}
Suppose that the assumptions of Theorem 2.5  are satisfied for $ V, F, f, f_n, n=1,2, \cdots $, where we 
require that $ V \in W_{1,2}((0, \delta))$ for some $\delta >0$ only if the  $f_n, n=1,2,\cdots $ are not 
all identically zero for n large enough. Moreover, assume that for each fixed $x,t$, $ (\frac{\partial}{\partial t}F)(x,t,z)$ is
$C^1$ in the real sense and that, 
\beq
\left|\left(\acute{\frac{\partial}{\partial t}F}\right)(x,t,z)\right|
\leq \, C_{R,N} \,\,\mathrm{for}\,\, x\in [0, \infty ),\, t\in I_N, \, |z|\leq R.
\label{3.45b}
\ene
Note that (\ref{3.45b}) implies (\ref{3.35f}). Then, the solution \, $u\in C\left([0,T_0], \3 \right), u(0)= \phi,\,  T_0 \leq T$, 
to the FNLSP (\ref{1.1}) 
depends continuously on the initial value and
on the boundary condition. In a precise way, let $u\in C\left([0, T_0], \3 \right)$ be the 
solution
to the FNLSP (\ref{1.1}) with $u(0)= \phi$. Let $ \phi_n \rightarrow \phi$ in $\3$  satisfy,
 $ \phi_n(0)= f_n(0)$ and
 assume that $f, f_n \in C^{3}$ and that
$ f_n(t) \rightarrow f(t)$ in 
$C^3([0,T_0])$. Moreover, if all the $f_n$ are not identically zero for $n$ large enough, suppose that $ \frac{\partial^2}{\partial t^2}F(0,t,z), 
\frac{\partial^2}{\partial t \partial z}F(0,t,z)$ and $ \acute{\acute{F}}(0,t,z)$ are 
continuous
in $z$, uniformly in $t$.
Then, for $n$ large enough the solution $u_n \in C\left([0,T_0], \3 \right) $ to the 
FNLSP (\ref{1.1}) with initial value
$\phi_n$ and boundary condition $f_n$ exists for $t \in [0,T_0]$ and $u_n \rightarrow u$ in 
$C\left([0,T_0], \3 \right)$.  
\end{theorem}

\noindent {\it Proof:}\,\, As in the proof of Theorem 2.3 it is enough to prove a local 
version for $T_0$ small enough. We define
$v_{0,n}$ and $P_n$ as in the proof of Theorem 2.3. As in the proof of Theorem 2.5 we prove that for $R$ large enough and $T_0$ small enough
all the $P_n$ send ${\cal N}_R$ (where we now require that $v(0)= \ v_{0,n} $) into 
$\overline{{\cal N}_R}$ and are contractions in the
norm of $B$ with a uniform contraction rate $\sigma < 1$ independent of $n$. Let $v_n$
be the unique fixed point. Then, $ u_n:= v_n+ r_n \in C\left([0,T_0], \3 \right)$
are the unique solutions to the FNLSP (\ref{1.1}) with initial value $\phi_n$ and boundary 
condition $f_n$.
Furthermore, 
$$
\|v_n -v \|_{B}=\left\|e^{-it H}[v_{0,n} -v_0]  + 
\frac{1}{i} \left(G {\cal F}_n(v_n)- G {\cal F}(v)\right)\right\|_{B} \leq \|v_{0,n}-v_0 \|_{L^2}+
 \sigma \left\| v_n - v \right\|_{B}+ 
$$
\beq
C \, T_0 \|f_n -f\|_{C^2([0,T_0])}.
\label{3.48}
\ene                      
In consequence,  $ v_n \rightarrow v$ in $B$. By taking the derivative of (\ref{3.9})
with respect to $t$ we obtain that,
\beq
i\frac{\partial}{\partial t}v = \u [H v_0 +{\cal F}( v_0)]+  G E,
\label{3.49}
\ene
where,
$$
E\left(v, \frac{\partial}{\partial t}v\right):=\acute{F}(x,t, v+r)\left(\frac{\partial}{\partial t}v+ 
\frac{\partial}{\partial t}r\right)+ 
\left(\frac{\partial}{\partial t}F\right)(x,t, v+r)
- i \frac{\partial^2}{\partial t^2}r +
$$
\beq 
\frac{\partial}{\partial t}r  \,V(x)- \frac{\partial^3}{\partial t \partial x^2}r. 
\label{3.50}
\ene

We similarly prove that $\frac{\partial}{\partial t}v_n$ satisfies,
\beq
i\frac{\partial}{\partial t}v_n= \u [H v_{0,n} +{\cal F}_n(v_{0,n})]+  G E_n,
\label{3.51}
\ene
with,
$$
E_n\left(v_n, \frac{\partial}{\partial t}v_n\right):=\acute{F}(x,t, v_n+r_n)\left(
\frac{\partial}{\partial t}v_n+ \frac{\partial}{\partial t}r_n \right)+ 
\left(\frac{\partial}{\partial t}F\right)(x,t, v_n+r_n)
- i\frac{\partial^2}{\partial t^2}r_n + 
$$
\beq
\frac{\partial}{\partial t}r_n  \, V(x)- \frac{\partial^3}{\partial t\partial x^2}r_n.  
\label{3.52}
\ene
By (\ref{3.49}) and (\ref{3.51}),
$$
i\frac{\partial}{\partial t}v_n(t)-i\frac{\partial}{\partial t}v(t)= \u [H v_{0,n} +{\cal F}_n
(v_{0,n})- H v_0 -{\cal F}(v_0)] 
+ G E_n\left(v_n, \frac{\partial}{\partial t}v_n \right)- 
G E_n\left(v_n, \frac{\partial}{\partial t}v\right)+
$$
\beq
 G E_n\left(v_n, \frac{\partial}{\partial t}v \right)- 
G E\left(v,\frac{\partial}{\partial t}v\right). 
\label{3.53}
\ene
Furthermore,
\beq
\left\|G \left[E_n\left(v_n,\frac{\partial}{\partial t}v_n \right)- 
E_n\left(v_n,\frac{\partial}{\partial t}v\right)\right] \right\|_{B}\leq C \,T_0 \left\|\frac{\partial}{\partial t}
v_n- \frac{\partial}{\partial t} v\right\|_{B}.
\label{3.54}
\ene

Hence, by (\ref{3.53}) and (\ref{3.54}) if  $ C\,  T_0 < 1/2   $,
$$
\lim_{n \rightarrow \infty}\left\|\frac{\partial}{\partial t} v_n- \frac{\partial}{\partial t} 
v \right\|_{B} \leq 2 \lim_{n \rightarrow \infty} \left[\left\| H v_{0,n} +{\cal F}_n
(v_{0,n})- H v_0 -{\cal F}(v_0) \right\|_{L^2}+  \left\|G \left[E_n  \left(v_n,
\frac{\partial}{\partial t}v \right)- \right.\right.\right.
$$
\beq
\left.\left.\left.
E\left(v, \frac{\partial}{\partial t}v \right)\right] \right\|_{B}\right]=0,
\label{3.55}
\ene
where we used that as $ v_n \rightarrow v$ in $L^2$ and $\| v_n \|_{\0}\leq C$, it follows by 
interpolation
\cite{rs} that $v_n \rightarrow v$ in $W_{s,2}, 0 < s <1$, and Sobolev's theorem.
Then, by the FNLSP (\ref{1.1})  $Hv_n \rightarrow Hv$ in $L^2$, and since $v_n \rightarrow v$ in $L^2$,  we have that 
$v_n \rightarrow v$ in the norm of $\1= D\left(\sqrt {H+M}\right)$. It follows that $v_n$ converges to $v$ in the norm of $\4$, and then  
$ u_n \rightarrow u$ in $  C\left([0,T_0], \3\right)$. In a standard way 
we extend the result of the theorem -step by step- to the original interval. For this purpose it 
is essential that the interval  of existence given by Theorem 2.5 depends only on the $\3$ norm of $\phi$.

\begin{remark}
{\rm
We prove as in Remark 2.4
that if  $I=[0, \infty)$  the solution  in $\3$ exits for all times unless it blows 
up in the $\3$ norm for some finite time, and that Theorem 2.5 
 implies  that the FNLSP (\ref{1.1}) has a unique solution 
$u \in C\left( I, \3 \right)$, with $u(0)= \phi \in \3$.}
\end{remark}

\bull

If the assumptions of Theorem 2.5 are satisfied, for any $\phi \in \3$ 
the FNLSP (\ref{1.1}) has a unique solution in $\2$ and a unique solution in $\3$
both with $ u(0)= \phi$. In the proposition below we prove that it is impossible that the
$\3$ solution blows up before the $\2$ solution does.
\begin{prop}(Regularity)
Suppose that the assumptions of Theorem 2.5 are satisfied.
Let $u \in C\left([0,T], \2\right)$ be a solution to the FNLSP (\ref{1.1})with
$u(0)= \phi \in \3$. Then $ u \in C\left([0,T], \3\right)$.
 \end{prop}

\noindent {\it Proof:}\,\, By Theorem 2.5 there is $ T_m \leq T $ such that $u(t)\in
  C\left([0,T_m], \3\right)$ and $\frac{\partial}{\partial t }u 
\in C\left([0,T_m], L^2\right)$. 
Furthermore, $\frac{\partial}{\partial t}v$ is a solution of the real-linear equation
(where $v$ is now fixed) (\ref{3.49}). Applying the contraction mapping theorem - step by step-
to this equation we prove that  $\frac{\partial}{\partial t}v \in C([0,T], L^2) $, and then,
it follows from  equation (\ref{3.7}) that $u = v +r \in C\left([0,T], \3\right)$.

\bull

We impose now further restrictions on $F$ that will allow us to derive an a-priori bound on the
 $\2$ norm of the solutions, and then, by Remark 2.4 that the solutions exist for all
times. We say that $F$ satisfies the sign condition if
\beq
\hbox{Im}\, \overline{z}F(x,t,z)=0, \, x, t \in \ER^+, z \in {\bf C},
\label{3.58}
\ene 
and we say that there is a hamiltonian structure if there is a function $h(x,t,z)$, such that 
for each fixed $x,t \in \ER^+$, $h$ is in 
$C^2\left({\bf C}, \ER\right)$ in the real sense, $h(x,t,0)=0$ and,
\beq
F(x,t,z)= 2 \frac{\partial}{\partial \overline{z}}h(x,t,z).
\label{3.59}
\ene
If Assumption A is satisfied we have that,
\beq
\left|h(x,t,z)\right|+ \left|\left(\frac{\partial}{\partial x}h\right)(x,t,z)\right|
\leq \, C_{R,N}\, |z|^2, \,\,\hbox{for}\,\, x\in [0, \infty),\, t\in I_N, \, |z|\leq R.
\label{3.59b}
\ene

Remark that as $h(x,t,0)=0$, equation (\ref{3.35f}) implies that $h(x,t,z)$ is differentiable in $t$, and that for 
each $R >0$ 
and each 
bounded
subset, $I_N$, of $I$,   
there is a constant $C_{R,N}$ such that,
\beq
\left|\left(\frac{\partial}{\partial t}h\right)(x,t,z)
\right| \leq \, C_{R,N}\, |z|^2,\,\hbox{for}\,\, x\in \ER^+,\, t\in I_N, \, |z|\leq R.
\label{3.60}
\ene 

Note that if (\ref{3.59}) is satisfied, then, (\ref{3.58}) is true if and only if $h$ depends 
only on $|z|$, i.e., if $h(x,t,z)=h(x,t,|z|)$ \cite{ka}.

Below we always assume that $\3 \subset W_{2,2}$.

 For any solution 
$u \in C\left(I, \3\right)$ to the
FNLSP (\ref{1.1}) the following  identities hold.
If (\ref{3.58}) is satisfied,

\beq
\frac{d}{dt}\|u(t) \|^2_{L^2}= 2 \hbox{Im}P(t)\overline{f(t)},
\label{3.61}
\ene
where we denote, $P(t):= (\frac{\partial}{\partial x}u)(0,t)$.
Observe that in the case where there is no external force, $f \equiv 0$, this is the conservation
of the $L^2$ norm.
Moreover, let $W(t)$ be the Hamiltonian,
\beq
W(t):= \frac{1}{2} \left\| \frac{\partial}{\partial x}u(t)\right\|^2_{L^2}+ \int_{\ER^+}
\left(\frac{1}
{2}V(x)|u(x,t)|^2+   h(x,t,u)\right)\, dx.
\label{3.62}
\ene
Then, if (\ref{3.59}) is true,
\beq
\frac{d}{d t} W(t)= - \hbox{Re}\acute{f}(t)\overline{P(t)} +\int_{\ER+} 
\left(\frac{\partial}{\partial t}h\right) (x,t,u )\,dx.
\label{3.63}
\ene
In the case where there is no external force and $h$ is independent of time  this 
identity is  the conservation of energy.
Furthermore, if (\ref{3.59}) is satisfied, 
\beq
 \frac{d}{d t}(u,u_x)= -i |P(t)|^2+ 2 i h(0,t,f(t)) - f(t)\,
\overline{\acute{f}(t)}
- 2 i \hbox{Re} \left(V u,u_x \right)+ 2 i \int_{\ER^+}\left(\frac{\partial}{\partial x}h\right)(x,t,u) \, dx.
\label{3.64}
\ene
Note that if $V$ is differentiable,
\beq
- 2 \hbox{Re} \left(V u,u_x \right)=  V(0)|f(0)|^2 + \int_{\ER^+} \acute{V}(x)|u(x,t)|^2 \, dx.
\label{3.64b}
\ene
The identity (\ref{3.64}) is analogous to the conservation of  momentum in the pure initial value
problem in $\ER$, c.f., \cite{gi}. Remark, however, that even in the case without external force
and with potential, $V$, and $\frac{\partial}{\partial x}h(x,t,u)$ both identically zero it is not a conservation law. This is to be expected 
because our problem is not translation invariant. 
The identities (\ref{3.61}), (\ref{3.63}) and(\ref{3.64}) where proven in the case
 $V \equiv 0$ and with  $F$ a 
single power, $F=\lambda |u|^{p-1} u$ in \cite{b1} and \cite{cb} (see also \cite{bs} for the multidimensional case)
for suitable smooth solutions. For the reader's convenience, we briefly give  below 
the details that show that  
the proof extends to our case, and that it  holds for solutions 
$u \in C\left(I, \3 \right)$. As $u(t) \in W_{2,2}, \lim_{x \rightarrow  \infty} u(x,t)=
 \lim_{x \rightarrow  \infty} \frac{\partial}{\partial x}u(x,t)= 0$. Then, by (\ref{1.1}),
(\ref{3.58}) and integrating by parts,                    
\beq
\frac{d}{d t}\| u(t)\|^2_{L^2}= 2 \hbox{Re}\left( \frac{\partial }{\partial t}u(t), u(t)\right)
= 2 \hbox{Im}\left(-\frac{d^2}{dx^2}u(t), u(t)\right)= 2 \hbox{Im}\, \overline{u(0,t)}\,
\frac{\partial}{\partial x}u(0,t),
\label{3.65}
\ene       
and (\ref{3.61}) holds. 
Moreover, denoting $u_x:= \frac{\partial}{\partial x}u$ and $ u_{xx} = 
\frac{\partial^2}{\partial x^2}u$,  we have that,
$$
\frac{\partial}{\partial t} \frac{1}{2}(u_x(t), u_x(t))= \lim_{\delta \rightarrow 0}\, 
\frac{1}{2 \delta}
\left[(u_x(t+\delta)- u_x(t), u_{x}(t+\delta))+( u_{x}(t), u_x(t+\delta)- u_x(t))\right]=
$$
\beq
-  \hbox{Re}\left(\frac{\partial}{\partial t}u(t),u_{xx}(t)\right)- \hbox{Re}\acute{f}(t)\overline{P(t)},
\label{3.66}
\ene  
where we integrated by parts before taking the limit $ \delta \rightarrow 0$.

Hence, by (\ref{1.1}) and (\ref{3.59}),

$$ 
\frac{d}{d t} W(t)= \hbox{Re}\left(\frac{\partial}{\partial t}u(t), 
i\frac{\partial}{\partial t}u(t)\right)
- \hbox{Re}\acute{f}(t)\overline{P(t)} + \int_{\ER+} \left(\frac{\partial}{\partial t}h
\right) (x,t,u )\,dx 
= 
$$
\beq
- \hbox{Re}\acute{f}(t)\overline{P(t)} +\int_{\ER+} \left(\frac{\partial}{\partial t}h
\right) (x,t,u )
\,dx,
\label{3.67}
\ene
and (\ref{3.63}) holds. Finally, integrating by parts, and using (\ref{1.1}),
$$
\frac{d}{d t}(u,u_x)= \lim_{\delta \rightarrow 0}\, 
\frac{1}{ \delta}
\left[(u(t+\delta)- u(t), u_{x}(t+\delta))+(u(t), u_x(t+\delta)- u_x(t))\right]=
2 i \hbox{Im}\left(\frac{\partial}{\partial t }u(t), u_x\right)- f(t)\,\overline{\acute{f}(t)}=
$$
$$
-2i \hbox{Re} ( H u+ F(x,t,u), u_x)- f(t)\,\overline{\acute{f}(t)}=
-i |P(t)|^2+ 2i h(0,t,f(t)) - f(t)\,\overline{\acute{f}(t)}
$$
\beq
 -2i \hbox{Re}(V u,u_x)+2 i \int_{\ER^+}\left(\frac{\partial}{\partial x}h
\right)(x,t,u)\, dx,
\label{3.68}
\ene
and (\ref{3.64}) holds. For any function $f$ we denote by $f_+$ its 
positive part and by $f_-$ its negative part, i.e., $f =f_+ - f_-, f_{\pm} \geq 0$. Below we denote by $\acute{V}$ the derivative of $V$ in distribution sense.

\begin{theorem}
Suppose that the assumptions of Theorem 2.5 are satisfied with $I=[0,\infty)$,  that $ \3 \subset W_{2,2}$ and  that   $\acute{V}$ is a function with,       
\beq
(\acute{V})_+  \leq C V_1+ Q, \; \mathrm{and} \; (\acute{V})_ - \in L^1_{\mathrm{loc}}(\ER^+),
\label{3.68b}
\ene 
where   $ Q$ satisfies (\ref{2.1}). Furthermore, assume  
that (\ref{3.58}), and (\ref{3.59}) 
 hold, where for each fixed $x,t \in \ER^+$, $h$ is in 
$C^2\left({\bf C}, \ER\right)$, in the real sense, and $h(x,t,0)=0$.
Moreover, assume that for each bounded subset
$I_N$ of $I$ there is a constant $C_N$ such that,
\beq
\left(\frac{\partial}{\partial x}h\right)_+(x,t,z) \leq \, C_N\, |z|^2,\,\mathrm{for}\,\,
 x\in \ER^+,\, t\in I_N, \, z \in {\bf C},
\label{3.69}
\ene 
and that for some $ 1 < p \leq 3$, 
\beq
\left(\frac{\partial}{\partial t}h\right)_+(x,t,z)
 \leq \, C_N\, \left(|z|^2 + |z|^{p+1} \right),\,\mathrm{for}\,\, x\in \ER^+,\, t\in I_N, \, z \in {\bf C},
\label{3.70}
\ene 
and, 
\beq
 h(x,t, z) \geq  - C_N (|z|^2+ |z|^{p+1}),\, 
\mathrm{for}\,\, x\in \ER^+, t \in I_N,
z \in {\bf C}.
\label{3.71}
\ene  
Then, the solutions in $\2$   and in $\3$ to the FNLSP (\ref{1.1}) given, respectively,   by
Theorems 2.2 and 2.5  exist for all time $t \in [0,\infty )$. 
\end{theorem}

\noindent {\it Proof:} In view of Remark 2.4 and of Proposition 2.8 it is enough to 
prove that for any finite
time interval $[0, T)$, the solution $u \in C\left([0,T), \2 \right)$ 
remains bounded in the norm of $\2$, as $t \rightarrow T$.
Suppose first that the solution $u \in C\left([0,T), \3 \right)$. For $ 0 \leq  t_1\leq t  <  T$ we denote,
$a(t):= (\int_{t_1}^t |P(\tau)|^2\, d\tau )^{1/2}$. In the estimates below we designate by
$C_T$  any constant that depends only on $T$  and $f$, and by $C_{T,1}$ any constant that depends on  $T$, $f$, and on the norm  $\|u(t_1)\|_{\2}$.
We denote, 
\beq
b( \phi ) := \mathrm{max} [\|\acute{\phi} \|_{L^2},  \|q \phi \|_{L^2}].
\label{3.71b}
\ene
 By (\ref{3.61})
\beq
\|u(t)\|^2_{L^2} \leq \|u(t_1)\|^2_{L^2} + C_T \ a(t), t_1 \leq t < T.
\label{3.72}
\ene
We denote,
\beq
\alpha(t):= \sup_{t_1 \leq s \leq t} b( u(s)).
\label{3.72b}
\ene
Integrating (\ref{3.64}) from $t_1$ to $t$, using (\ref{2.2}),  (\ref{3.64b}), (\ref{3.68b}), (\ref{3.69}), (\ref{3.72}) and  the estimate 
$ |(u,u_x)|\leq 1/2 \|u\|_{L^2}^2 + 1/2 \|u_x\|_{L^2}^2$, we prove that,
\beq
a(t)\leq   C_T \alpha(t)+ C_{T,1} , t \in [t_1,T),
\label{3.73}
\ene
where we used that $ a(t)$ is a non-decreasing function. Integrating again (\ref{3.64}) from $t_1$ to $t$,  using now 
(\ref{2.2}), (\ref{3.64b}), (\ref{3.68b}), (\ref{3.69}), (\ref{3.72}), (\ref{3.73}) and as 
$|(u, u_x)|\leq \|u\|_{L^2} \, \|u_x\|_{L^2}$ we obtain that, 

\beq
a(t) \leq  \sqrt{\|u(t)\|_{L^2}\,\|u_x(t)\|_{L^2}} + C_{T,1} + C_T (T-t_1)  \sqrt{\alpha(t)},     t_1 \leq t < T.
\label{3.74}
\ene
We denote, $ g(t_1,t):=(\int_{t_1}^{t} |f(\tau)|^2\, d\tau )^{1/2}$. 
Now we integrate (\ref{3.61}) from $t_1$ to $t$, and using (\ref{3.74})  we prove that,
\beq
\|u(t)\|_{L^2}^2\leq \|u(t_1)\|_{L^2}^2 + 2 g(t_1,t)\left[\sqrt{\|	u(t)\|_{L^2}\,\|u_x(t)\|_{L^2}}+ C_{T,1}+
C_T  (T-t_1)  \sqrt{ \alpha(t)} \right] , t_1 \leq t <T.
\label{3.76}
\ene
By (\ref{3.76}) for some constant $C$,
\beq
\left\|u(t)\right\|_{L^2}^2 \leq  C \left[ \left\|u(t_1)\right\|_{L^2}^{8/3} + 2^{4/3}\, (g(t_1,t))^{4/3} \, 
\left( \|u_x(t)\|^{2/3}_{L^2}+ C_T (T-t_1)^{4/3}  \alpha(t)^{2/3}+ C_{T,1} \right) \right] +1,
  \,  t_1 \leq t <  T. 
\label{3.77}
\ene
Here we consider first the case, $\|u(t)\|_{L^2} \leq 1$, where the estimate is trivial, and then the case  $\|u(t)\|_{L^2} \geq 1$. 
For any $ 1 \leq p < 5 $ there is a constant $C$ such that for any $u \in W_{1,2}$ and any $\epsilon > 0$,     
\beq
\| u \|_{L^{p+1}}^{p+1}\leq C \epsilon \| u_x\|_{L^2}^2 + \frac{C}{\epsilon^{(\nu-1)/2}}
\| u\|_{L^2}^{2\nu},
\label{3.78}
\ene
where $\nu:= 1+ \frac{2(p-1)}{5-p}$. 
We give the proof of (\ref{3.78}) below.   
Integrating (\ref{3.63}) from $t_1$ to $t$ , and by (\ref{2.2}) with $\epsilon=1/4$,  (\ref{3.70}), (\ref{3.71}), (\ref{3.72}), (\ref{3.73}),  (\ref{3.77}), (\ref{3.78}) 
and as $ \nu \leq 3$,
\newpage
$$
\frac{1}{4} b( u(t))^2 
\leq W(t_1)+ C_{T,1} + C_T  \alpha(t)  + C_T (T-t_1)\left[ \epsilon  
\alpha(t)^2 + \frac{1}{\epsilon^{(\nu-1)/2}} \left(C_{T,1} +  \right.\right.
$$
\beq 
\left.\left.
  g(t_1,t)^{4\nu /3} \alpha(t)^2\right)\right], t_1 \leq t < T. 
\label{3.79}
\ene
Pick any $ \epsilon $ and   $\Delta$ such that, 
$ C_T T [\epsilon + \frac{1 }{\epsilon^{(\nu-1)/2}}   \|f\|_{L^{\infty}([0,T])}^{4\nu/3} \Delta^{2\nu/3}\leq 1/8.$ Then, 
by (\ref{3.79}),
\beq
\alpha(t)^2 \leq 8\left\{ W(t_1)+ C_{T,1} + C_T  \alpha(t) + \frac{C_{T,1}}{\epsilon^{(\nu-1)/2}}\right\}, t_1 \leq t <   \min [t_1+\Delta, T]. 
\label{3.80}
\ene
As by Theorem 2.3 and Proposition 2.8 we can approximate solutions in $\2$ by solutions in $\3$,
equations (\ref{3.76}), (\ref{3.77}), (\ref{3.79}) and  (\ref{3.80}) hold  also if $u \in
 C\left([0,T), \2\right)$.

Suppose now that we are given a solution $ u \in C\left([0,T_m), \2\right)$  to the FNLSP (\ref{1.1})
where  $T_m$ is the maximal time of existence. Then, we must have $T_m = \infty$, because
if $T_m < \infty$  we can take $t_1= T_m - \Delta $, and then by (\ref{3.72}), (\ref{3.73}) and  (\ref{3.80}),
\beq
\| u(t)\|_{\2} \leq C_{T,1}, \; \mathrm{for}\; T_m- \Delta  \leq t < T_m ,
\label{3.81}
\ene   
 and  by Remark 2.4 we can continue $u(t)$ to $t > T_m$, in contradiction with the 
definition of $T_m$. 

We now prove (\ref{3.78}). By the Sobolev-Gagliardo-Nirenberg inequality  \cite{fr}
\beq
\|u\|_{p+1}^{p+1} \leq  C \| u_x\|_{L^2}^{a (p+1)}\, \|u\|_{L^2}^{(p+1)(1-a)}, \, a:= \frac{1}{2}- \frac{1}{p+1}.
\label{3.82}
\ene
 Inequality (\ref{3.82}) is stated in \cite{fr} for $u \in C^{\infty}_0(\ER)$, but by continuity it applies to $u \in W_{1,2}(\ER)$ and
 extending $u \in W_{1,2}$ as an even function in $W_{1,2}(\ER)$ it also holds for $u \in W_{1,2}$.
 Denote, $ k:= \frac{5-p}{4}$. Then, by (\ref{3.82}),

\beq
\|u\|_{p+1}^{p+1} \leq  C \| u_x\|_{L^2}^{2(1-k) }\, \|u\|_{L^2}^{ 2 \nu k}
 \leq C \epsilon \| u_x\|_{L^2}^2 + \frac{C}{\epsilon^{(\nu-1)/2}}
\| u\|_{L^2}^{2\nu},
 \label{3.83}
\ene
where we used the inequality, $a^{1-k}\, b^k \leq \epsilon a + \frac{1}{\epsilon^{(1/k -1)} } b,  a,b \geq 0,  \epsilon > 0, 0 < k \leq 1.$

\newpage

\begin{remark}
{\rm
In the case where $f \equiv 0$ we prove that the solutions are global under weaker assumptions
because we do not need to use identity (\ref{3.64}).
Suppose that assumptions A and B are satisfied with $I=[0,\infty)$, and with $f \equiv 0$,
that (\ref{3.35f}), (\ref{3.58}), and (\ref{3.59}) hold, where for each fixed $x,t \in \ER^+$, $h$ is in 
$C^2\left({\bf C}, \ER\right)$ in the real sense, and $h(x,t,0)=0$.
Moreover, assume that $\4 \subset W_{2,2}$, that   (\ref{3.70}) and (\ref{3.71}) hold with $ 1 < p < 5$. Then, the conclusions of Theorem
2.9 are true.
The proof is much simpler now because by (\ref{3.61}) $\|u(t) \|\leq C$, and then 
by (\ref{3.63}), (\ref{3.70}), (\ref{3.71})  and (\ref{3.78}), $b(u(t)) \leq C$. We complete the proof as in
Theorem 2.9.} 
\end{remark}
\begin{remark}
{\rm
Recall that  our results in local solutions given in Theorems 2.2, 2.3, 2.5 and 2.6, 
in Remarks 2.4, 2.7 and in Proposition 2.8  hold without any restriction on the grow of $V_1$ 
at infinity. Our results in global solutions given in Theorem 2.9 and in Remark 2.10  require that $\3 \subset W_{2,2}$. We give now a sufficient condition for this to hold.
We denote by $\mathit{Lip}$ the set of all continuous and bounded functions, $f$, defined on $[0,\infty)$ that are globally Lipschitz, i.e. such that,
\beq
\mathit{Lip}(f):= \sup_{x,y\in [0,\infty), x \neq y } \frac{\left| f(x)-f(y) \right|}{|x-y|} < \infty.
\label{3.84}
\ene
Note that if $ f \in \mathit{Lip}$ then $f$ is differentiable for a.e. $x$ with $\acute{f}\in L^{\infty}$ and $\mathit{Lip}(f):= \|\acute{f}\|_{L^{\infty}}$. Suppose as above that
$V=V_1+V_2, V_j \in L^1_{\mathrm loc}(\ER^+), j=1,2, V_1 \geq 0$ and $V_2$ satisfies (\ref{2.3}). Remark that, eventually adding $1$ to $V_1$ and substracting 
it from $V_2$, we can
assume that $V_1 \geq 1$. Suppose that $g:= (V_1)^{-1/2} \in \l$. For $c >0$ denote $g_c:= (V_1 +c)^{-1/2}$. Observe that, $ \acute{g}_c:= (1+c V_1^{-1})^{-3/2}\, \acute{g}$.
Then, $ g_c \in \l$ and $\l(g_c)$ decreases monotonically as $c\rightarrow \infty$. It follows from Theorem 7.1 of \cite{ka1} (this paper considers the case in the whole line, 
but the proof in our case is the same) that if $\lim_{c \rightarrow \infty} \l{g_c} < 1$, then $H_1:=H_0+V_1$ is selfadjoint in the domain,  $D(H_1)=W_{2,2}\cap\0 \cap
D(V_1)$ . Hence, by (\ref{2.4}) and Kato-Rellich's theorem,  $D(H)= W_{2,2}\cap\0 \cap D(V_1)$. Assume moreover, that 
$V_1 \in L^2_{\mathrm{loc}}([0,\infty))$. Let us take any $h \in C^{\infty}_0([0, \infty))$, satisfying 
$ h(x)=1, 0 \leq x \leq 1 $. We decompose any 
$ \phi \in \3$ as $\phi = \phi_1+\phi_2$, with $\phi_1:= \phi -\phi(0) h, \phi_2:= \phi(0) h$. Then, under the assumptions above, $\phi \in \3 \leftrightarrow \phi_1 \in D(H)$,
and it follows that in this case $\3 = W_{2,2} \cap D(V_1)$. Note that  $V_1$ can be any positive polinomial, $p(x)$, or $\mathrm{exp}p(x), \mathrm{exp}\,\mathrm{exp}p(x), \cdots$. 
 Moreover, (\ref{3.68b})  is satisfied, for example, if $V_1$ any positive polynomial, or $V_1=e^ x$,  and $(\acute{V_2})_+$ fulfills (\ref{2.1}). Finally, note that in
 the case of Remark 2.10 where the force is identically zero we do not need that $V_1 \in L^2_{\mathrm{loc}}([0,\infty))$. In this case we can admit, for example,
 $V_1= \frac{1}{x^k}, 0< x
\leq 1, V_1= 1, x \geq 1, k  \geq 2$ (for the case k=2 see Example 7.4 of \cite{ka1}).
}
\end{remark}

\end{document}